\documentclass[12pt]{amsart}
\usepackage{amssymb}
\newtheorem{theorem}{Theorem}[section]
\newtheorem{proposition}[theorem]{Proposition}
\newtheorem{coro}[theorem]{Corollary}
\newtheorem{lemma}[theorem]{Lemma}

\newcommand\N{\mathbb{N}}
\newcommand\R{\mathbb{R}}

\begin{document}

\baselineskip=17pt

\title[Riesz meets Sobolev]
      {Riesz meets Sobolev}

\author{Thierry Coulhon}
\author{Adam Sikora}
\address{Thierry Coulhon, D\'epartement de Math\'ematiques, Universit\'e
de
Cergy-Pontoise, Site de Saint-Martin, 2, rue Adolphe Chauvin, F 95302
Cergy-Pontoise Cedex, France}
\email{coulhon@u-cergy.fr}
\address{Adam Sikora, Department  of Mathematics, Macquarie University, NSW 2109, Australia
}
\email{adam.sikora@anu.edu.au}
\date{\today}
\subjclass[2010]{58J35, 42B20, 46E35}
\keywords{Heat kernels, Riesz transform, Sobolev type inequalities}

\begin{abstract}
We show that the $L^p$ boundedness, $p>2$, of the Riesz transform on a complete non-compact Riemannian manifold with upper and lower Gaussian heat kernel estimates is equivalent to a certain form of Sobolev inequality.
We also characterize in such terms the  heat kernel gradient upper estimate  on manifolds with polynomial growth.
\end{abstract}

\maketitle
\setcounter{section}{0}
\renewcommand{\theenumi}{\alph{enumi}}
\renewcommand{\labelenumi}{\textrm{(\theenumi)}}
\numberwithin{equation}{section}

{\it In memoriam Nick Dungey and Andrzej Hulanicki.}

\tableofcontents

\section{Introduction}

The present paper may be considered as a companion paper to \cite{ACDH}, which gave criteria for the $L^p$ boundedness, for $p>2$, of the Riesz transform on non-compact Riemannian manifolds. Here we reformulate these criteria in terms of certain Sobolev inequalities. That is, we deduce some  $L^p$ to $L^p$ estimates from  suitable $L^q$ to $L^p$ estimates, for $q<p$.

Let  $M$ be a complete, connected, non-compact Riemannian manifold. The methods of this paper remain valid for other types of
spaces endowed with a gradient, a metric which is compatible with this gradient, a measure, and finally an operator associated
with the Dirichlet form constructed from the gradient and the measure. An interesting example is a Lie group endowed with a family
of left-invariant H\"ormander vector fields. We leave  the details of such extensions to the reader.

Let $d$ be the geodesic distance on $M$; denote by $B(x,r)$  the open ball with respect to $d$ with center $x\in M$ and radius $r>0$.

Denote by  $\mu$ the Riemannian measure,   by $L^p(M,\mu)$, $1\le p\le \infty$, the corresponding $L^p$ spaces, and let
$V(x,r)=\mu(B(x,r))$.

Let  $\Delta$ be the  (non-negative) Laplace-Beltrami operator. One could consider another measure  ${\tilde\mu}$ with positive smooth non-zero density with respect to $\mu$,
and the associated operator $\Delta_{\tilde\mu}$, formally given
by $$(\Delta_{\tilde\mu}f,f)=\int_M|\nabla f|^2\,d\tilde\mu.$$ Again, for simplicity, we stick to the standard case.

Let $\nabla$ be the Riemannian gradient.
We can now define formally the Riesz transform operator $\nabla\Delta^{-1/2}$.

Let $p\in (1,\infty)$.
The boundedness of the Riesz transform on $L^p(M,\mu)$ reads
  $$
  \|\left|\nabla f\right|\|_p
\le  C_p \|\Delta^{1/2}f\|_p,\ \forall\,f\in {\bf
C}^\infty_0(M),\leqno{(R_p)}
  $$
and if the reverse inequality $(RR_p)$ also holds, one has
 $$
  \|\left|\nabla f\right|\|_p
\simeq  \|\Delta^{1/2}f\|_p,\ \forall\,f\in {\bf
C}^\infty_0(M).\leqno{(E_p)}
   $$

One says that $M$ satisfies the volume doubling property if there exists 
$C$ such that
$$V(x,2r)\le CV(x,r), \
\forall\,r>0,\,x\in M,\leqno{(D)}$$
more precisely if 
there exist $\nu, C_\nu>0$ such that
$$
\frac{V(x,r)}{V(x,s)}\le C_\nu\left(\frac{r}{s}\right)^\nu, \
\forall\,r\ge s>0,\,x\in M.
\leqno{(D_\nu)}
$$

The heat semigroup is the family of operators $(\exp(-t\Delta))_{t>0}$ acting  on $L^2(M,\mu)$, it has a positive and smooth  kernel $p_t(x,y)$
called the heat kernel.
In the sequel, we shall consider the following standard heat kernel estimates for manifolds with doubling :
the on-diagonal upper estimate,
$$
p_t(x,x) \le \frac{ C }{V(x,\sqrt{t})}\leqno{(DU\!E)}
$$
for some $C>0$,  all $x \in M$ and $t>0$,
the full Gaussian upper estimate,
$$
 p_t(x,y) \le \frac{ C }{V(x,\sqrt{t})}\exp\left(-c\frac{d^2(x,y)}{t}\right)\leqno{(U\!E)}
$$
for some $C,c>0$, all $x,y \in M$ and $t>0$,
 the upper and lower Gaussian estimates,
$$
\frac{ c }{V(x,\sqrt{t})}\exp\left(-\frac{d^2(x,y)}{ct}\right)\le p_t(x,y) \le \frac{ C }{V(x,\sqrt{t})}\exp\left(-\frac{d^2(x,y)}{Ct}\right)\leqno{(LY)}
$$
for some $C,c>0$, all $x,y \in M$ and $t>0$, and finally the gradient upper estimate
$$
|\nabla_x p_t(x,y)| \le \frac{C}{\sqrt{t}V(y,\sqrt{t})}\leqno{(G)}
$$
for all $x,y \in M$, $t>0$. It is known that, under $(D)$ and $(DU\!E)$, $(G)$ self-improves into $$
|\nabla_x p_t(x,y)| \le \frac{C}{\sqrt{t}V(y,\sqrt{t})}\exp\left(-\frac{d^2(x,y)}{Ct}\right)
$$
for all $x,y \in M$, $t>0$, see \cite{D} and also \cite[Section 4.4]{CS}). It will follow from Proposition \ref{gr} below that assumption $(DU\!E)$ is not needed here.

Recall that $(DU\!E)$  plus $(D)$  implies $(U\!E)$ (\cite[Theorem 1.1]{G2}, see also \cite[Corollary 4.6]{CS}), and $(U\!E)$ plus $(G)$ and $(D)$ implies $(LY)$
(\cite{LY}). We shall see   in Proposition \ref{gr} below that $(G)$ implies $(DU\!E)$, therefore $(G)$ plus $(D)$ implies $(LY)$. Conversely, $(LY)$ implies $(D)$ (see for instance \cite[p.161]{Sal1}). 

\bigskip

The following is one of the two main results  of \cite{ACDH} (Theorem~1.4 of that paper).

\begin{theorem}\label{maincor} Let $M$ be a complete non-compact Riemannian manifold satisfying $(D)$, $(DU\!E)$, and $(G)$. Then
  the equivalence $(E_p)$ holds for
$1<p<\infty$.
\end{theorem} 

Taking into account Proposition \ref{gr} below, one can skip condition $(DU\!E)$, and formulate this result in the following simpler way.

\begin{theorem}\label{maincorCS} Let $M$ be a complete non-compact Riemannian manifold satisfying $(D)$ and $(G)$. Then
  the equivalence $(E_p)$ holds for
$1<p<\infty$.
\end{theorem} 

Let us now introduce an $L^p$ version of $(G)$, namely
$$
\||\nabla e^{- t\Delta}|\|_{p\to p} \le  \frac{C_p}{\sqrt t},\ \forall\, t>0.\leqno{(G_p)}
$$

The other main result of  \cite{ACDH} is the following (Theorem~1.3 and Proposition~1.10 of that paper).

\begin{theorem}\label{nsc} Let $M$ be a complete non-compact Riemannian manifold satisfying $(D)$
and
 $(LY)$. Let $p_0 \in (2,\infty]$. The following assertions are
equivalent:
\begin{enumerate}
\item $(R_p)$ holds for all $p\in (2,p_0)$. 
\item
  $(G_p)$ holds for all $p \in (2,p_0)$. 
\item
 For all  $p \in (2,p_0)$, there exists $C_p$ such that  
$$
 \||\nabla p_{t}(.,y)|\|_{p}\le    \frac{C_p }
{\sqrt{t}\left[V(y,\sqrt{t})\right]^{1-\frac{1}{p}}},\ \forall\,t>0,\,y\in M.
$$
\end{enumerate}
\end{theorem}

According to Proposition \ref{prop12} below, we will be able
to add another equivalent condition in the above list, namely \begin{enumerate}
\item[(d)] {\it ~For all  }$p \in (2,p_0), $ {\it~ there exists }$ C_p$ {\it ~such that~}
$$ \|\left|\nabla f\right|\|_p^2
\le C_p  \|f\|_p\|\Delta f\|_p,\ \forall\,f\in {\bf
C}^\infty_0(M).$$
\end{enumerate}
\
The two above results  are the cornerstones of the present paper. Our main results are  Theorems \ref{lem4}, \ref{lem3} and Corollary \ref{E} below. In Theorem \ref{lem4},   using Theorem \ref{nsc}, we give a  necessary and sufficient condition for $(R_p)$  to hold for $p$ in an interval above $2$ on   manifolds with polynomial volume growth satisfying $(D)$
and
 $(LY)$, in terms of an $L^p-L^q$  Sobolev type inequality with a gradient in the left-hand side. In Theorem \ref{lem3},  we give a necessary and sufficient condition  for $(G)$   on   manifolds with polynomial volume growth  satisfying  a mild local condition, in terms of a multiplicative $L^\infty$ Sobolev type inequality, with a gradient in the left-hand side.  In Corollary \ref{E} we deduce from  Theorem \ref{lem4} and Theorem \ref{maincor}   that this $L^\infty$ Sobolev inequality alone implies $(E_p)$ on manifolds with polynomial growth and the above local condition.

Here is the plan we will follow. In section \ref{ghk}, we prove  that $(G)$ implies $(DU\!E)$, together with a similar statement for some related kernels. In section \ref{dv}, we give a first version of our results for manifolds with doubling and a polynomial volume upper bound. In section \ref{pg}, we assume full  polynomial growth and obtain more complete results. Finally, in section \ref{ap},
we give applications of our methods to  second order elliptic operators in $\R^n$. 

\section{Gradient estimates imply heat kernel bounds}\label{ghk}

Note that the following result does not require assumption $(D)$.

\begin{proposition}\label{gr}
$(G)$ implies $(DU\!E)$.
\end{proposition}

\begin{proof}
For $x\in M$, $t>0$, define  $$
K=K(x,t)= \frac{V(x,\sqrt{t})p_t(x,x)}{2}.
$$
We claim that
$$
p_t(y,x) \ge \frac{K}{V(x,\sqrt{t})}
$$
for all $y\in B(x,\frac{K\sqrt{t}}{C})$.
 Indeed, according to $(G)$ and the mean value theorem, for such $y$,
 $$
 |p_t(y,x)-p_t(x,x)| \le \frac{Cd(y,x)}{\sqrt{t}V(x,\sqrt{t})} \le
\frac{C}{\sqrt{t}V(x,\sqrt{t})}\frac{K\sqrt{t}}{C}=\frac{K}{V(x,\sqrt{t})}.
 $$
 Thus, given the definition of $K$,
 $$
 p_t(y,x)\ge p_t(x,x) -\frac{K}{V(x,\sqrt{t})} =
\frac{K}{V(x,\sqrt{t})},
 $$
 hence the claim.
 Now
 \begin{eqnarray*}
 1&\ge& \int_M p_t(y,x)\,d\mu(y)   \ge \int_{B(x,\frac{K\sqrt{t}}{C})}
p_t(y,x)\,d\mu(y) \\&\ge& \int_{B(x,\frac{K\sqrt{t}}{C})}
\frac{K\,d\mu(y)}{V(x,\sqrt{t})} =
  \frac{KV(x,\frac{K\sqrt{t}}{C})}{V(x,\sqrt{t})}.
 \end{eqnarray*}
 If  $K\ge C$, this means that $K\le 1$.
 Hence $$
 K \le \max \left(C,1\right),
 $$
and 
$$
p_t(x,x) \le \frac{ 2\max
\left(C, 1\right) }{V(x,\sqrt{t})}, \ \forall\,t>0,\,x\in M.
$$
\end{proof}

As we noticed in the introduction, the following is a consequence of Proposition \ref{gr} together with known results.

\begin{coro}\label{grd}
Assume that $M$ satisfies $(G)$ and $(D)$.
Then $M$ satisfies
$(LY)$.
\end{coro}

\bigskip

It may be of interest to notice that the assumption in Proposition \ref{gr} can be replaced by a gradient estimate of some other kernels.
Namely, for $a>0$,  denote by $r^a_t(x,y)$ the  (positive) kernel of the operator
$$R^a_t=(I+t\Delta)^{-a}=\frac{1}{\Gamma(a)}\int_0^{+\infty} s^{a-1}\exp(-s(I+t\Delta))\,ds.$$
Similarly, for $0<a<1$,  denote by $p^a_t(x,y)$  the kernel of the
operator $P^a_t=\exp(-(t\Delta)^a)$. In the following statement, we assume doubling only for simplicity,
otherwise one has to include an additional constant in the outcome. 

\begin{proposition}\label{gr1} Assume $(D)$.
Suppose that $q^a_t=r^a_t$ for some $a>0$ or $q^a_t=p_t^a$ for some
$0<a<1$.
Next assume that $M$ satisfies the gradient upper estimate
$$
|\nabla_x q^a_t(x,y)| \le \frac{C}{\sqrt{t}V(y,\sqrt{t})}\leqno{(G^a)}
$$
for all $x,y \in M$, $t>0$.
Then $M$ satisfies
$(DU\!E)$.
\end{proposition}

\begin{proof} First note that in all cases
$$\int_Mq_t^a(x,y)\,d\mu(y)\le 1,\ \forall\,x\in M.$$
Fix $x\in M$, $t>0$. Define  $$
K=K(x,t)= \frac{V(x,\sqrt{t})q^a_t(x,x)}{2}.
$$
Exactly as in the proof of Proposition \ref{gr}, one shows that
 $$
 K \le \max \left(C,1\right).
 $$
 To finish the proof of Proposition~\ref{gr1}, note that, for $a>0$,
 $$
 p_t(x,x)=\|p_{t/2}(\cdot, x)\|_2 \le C
\|r^{a/2}_{t/2}(\cdot, x)\|_2=Cr^a_{t/2}(x,x).
 $$
 Indeed, write
$$\exp(-t\Delta)= {\exp(-t\Delta)(I+t\Delta)^a}(I+t\Delta)^{-a} ,$$ 
so that $$p_t(\cdot, x)=  {\exp(-t\Delta)(I+t\Delta)^a}  r_t^a(\cdot, x),$$
and since by spectral theory the operator
$\exp(-t\Delta)(I+t\Delta)^a$ is uniformly bounded on $L^2(M,\mu)$, the claim is proved.

Similarly, for $0<a<1$, writing $$\exp(-t\Delta)= \exp(-t\Delta)\exp(t\Delta)^a\exp(-(t\Delta)^a),$$
one sees that 
  $$
 p_t(x,x)=\|p_{t/2}(\cdot, x)\|_2 \le C
\|p^{a}_{t/2}(\cdot, x)\|_2=Cp^a_{2^{\frac{1-a}{a}}t}(x,x).
 $$
 
\end{proof}

\section{Doubling volume}\label{dv}

Recall that $(LY)$ implies $(D)$. Thus the doubling volume assumption will be implicit in the first two statements of this section.

\begin{theorem}\label{lemd3}
Let $M$ satisfy  $(LY)$.
Let $\nu>0$ be such that 
\begin{equation}
V(x,r)\le Cr^\nu, \forall\,r>0,\,x\in M.\label{volnu}
\end{equation}
Let $p_0\in (2,\infty]$.
Assume
  \begin{equation}\label{sobpd}
  \|\left|\nabla f\right|\|_p
\le  C_p
\|\Delta^{\frac{\nu}{2}\left(\frac{1}{q}-\frac{1}{p}\right)+\frac{1}{2}}f\|_q,\
\forall\,f\in {\bf C}^\infty_0(M),
    \end{equation}
for  all $p\in (2,p_0)$ and some $1<q<p$.
Then $(R_p)$ holds for  all $p\in (2,p_0)$.
\end{theorem}

\begin{proof}
Let $p$ be such that $2<p<p_0$ and $q\in (1,p)$ such that   \eqref{sobpd} holds. Taking $f=p_{2t}(.,y)=\exp(-t\Delta)p_t(.,y)$, $t>0$, $y\in M$, in \eqref{sobpd}, one obtains
$$
 \||\nabla p_{2t}(.,y)|\|_{p}\le C_{p}
\|\Delta^{\frac{\nu}{2}\left(\frac{1}{q}-\frac{1}{p}\right)+\frac{1}{2}}{\exp(-t
\Delta)}p_t(.,y)\|_q,
$$
 hence, by analyticity of the heat semigroup on $L^q(M,\mu)$,
$$
 \||\nabla p_{2t}(.,y)|\|_{p}\le    C
t^{-\frac{\nu}{2}\left(\frac{1}{q}-\frac{1}{p}\right)-
\frac{1}{2}}\|p_t(.,y)\|_q ,\ \forall\,t>0,\,y\in M.
$$
On the other hand,
$(U\!E)$ yields
$$\|p_t(.,y)\|_q\le \frac{C}{\left[V(y,\sqrt{t})\right]^{1-\frac{1}{q}}},\
\forall\,t>0,\,y\in M.$$
Hence
$$
 \||\nabla p_{2t}(.,y)|\|_{p}\le    \frac{C t^{-\frac{1}{2}}}
{\left[V(y,\sqrt{t})\right]^{1-\frac{1}{p}}}\left[t^{-\frac{\nu}{2}} V(y,\sqrt{t})
\right]^{\frac{1}{q}-\frac{1}{p}},\,\ \forall\,t>0,\,y\in M,
$$
and,  according to \eqref{volnu},
the quantity $t^{-\frac{\nu}{2}} V(y,\sqrt{t})$ is bounded from above, therefore
\begin{equation}\label{glgl}
 \||\nabla p_{2t}(.,y)|\|_{p}\le    \frac{C' }
{\sqrt{t}\left[V(y,\sqrt{t})\right]^{1-\frac{1}{p}}},\,\ \forall\,t>0,\,y\in M.
\end{equation}
One concludes by applying \cite{ACDH},  namely Theorem \ref{nsc} above.

\end{proof}

{\bf Remarks:}

-Remember that it follows from \cite{CD} that, under the assumptions of Theorems \ref{lemd3}, $(R_p)$ also holds
for $p\in (1,2]$. As a consequence, $(RR_p)$ also holds, therefore assumption \eqref{sobpd}
implies  
$$  \|\Delta^{1/2} f\|_p
\le  C_p
\|\Delta^{\frac{\nu}{2}\left(\frac{1}{q}-\frac{1}{p}\right)+\frac{1}{2}}f\|_q,\
\forall\,f\in {\bf C}^\infty_0(M),
$$
hence,   by making the change of functions $\Delta^{1/2}f\to f$, the Sobolev inequality
$$  \|f\|_p
\le  C_p
\|\Delta^{\frac{\nu}{2}\left(\frac{1}{q}-\frac{1}{p}\right)}f\|_q,\
\forall\,f\in {\bf C}^\infty_0(M).
$$
It follows that 
$$V(x,r)\ge cr^\nu, \forall\,r>0,\,x\in M$$ (see \cite{Ca}).
Thus, in fact, under the assumptions of Theorem \ref{lemd3},  the volume growth of $M$ has to be polynomial of exponent $\nu$ (in particular, $\nu$ has to coincide with the topological dimension of $M$). However, the fact that we do not use explicitly polynomial growth in the proof will allow us
below some true excursions in the doubling volume realm.

-An equivalent formulation of \eqref{sobpd} is
$$
  \|\left|\nabla f\right|\|_p
\le  C_p
\|\Delta^{\alpha}f\|_q,\
\forall\,f\in {\bf C}^\infty_0(M),
$$
for  all $p\in (2,p_0)$ and some $\alpha>1/2$, with $q=\frac{1}{\frac{1}{p}+\frac{2}{\nu}\left(\alpha-\frac{1}{2}\right)}$.
In particular,  $\alpha=1$ and $q=\frac{p\nu}{\nu+p}$ is a valid choice. See Section \ref{ap} below.

-When $p_0<\infty$, if  one assumes
\begin{equation}
 \|\left|\nabla f\right|\|_{p_0}
\le  C
\|\Delta^{\frac{\nu}{2}\left(\frac{1}{q}-\frac{1}{{p_0}}\right)+\frac{1}{2}}f\|_q,\
\forall\,f\in {\bf C}^\infty_0(M),\label{rm}
\end{equation}
instead of \eqref{sobpd}, one still obtains the same conclusion by interpolation.

 -One can also replace \eqref{sobpd} by the following weaker inequality
$$
  \|\left|\nabla f\right|\|_p
\le  C_p
\|\Delta^{\alpha/2}f\|_{q_1}^{\theta}\|f\|_{q_2}^{1-\theta},\
\forall\,f\in {\bf C}^\infty_0(M),
$$
where $0<\theta <1$, $1\le q_1,q_2\le \infty$, $\frac{1}{p}<\frac{\theta}{q_1}+\frac{1-\theta}{q_2}<1$
and $$\alpha\theta=\nu\left(\frac{\theta}{q_1}+\frac{1-\theta}{q_2}-\frac{1}{p}\right)+1.$$ Here also, one can take $p=p_0$.

\bigskip

In the next statement, we shall relax the volume upper bound  assumption for small radii. This can be useful in situations where the volume growth is polynomial, but with different exponents for small and large radii, say for instance the Heisenberg group endowed with a group invariant Riemannian metric.

We shall say that the  local Riesz inequality $(R_p)_{loc}$ holds on $M$ if
 $$
  \|\left|\nabla f\right|\|_p
\le  C_p \left(\|\Delta^{1/2}f\|_p+\|f\|_p\right),\ \forall\,f\in {\bf
C}^\infty_0(M).  $$

This is the case for instance if $M$ has Ricci curvature bounded from below (see \cite{B}).

\begin{theorem}\label{lemd4}
Let $M$ satisfy  $(LY)$ and $(R_p)_{loc}$.
Let $\nu>0$ be such that 
\begin{equation}
V(x,r)\le Cr^\nu, \forall\,r\ge 1,\,x\in M.\label{volnu1}
\end{equation}
Let $p_0\in (2,\infty]$.
Assume
\eqref{sobpd}
for  all $p\in (2,p_0)$ and some $1<q<p$.
Then $(R_p)$ holds for  all $p\in (2,p_0)$.
\end{theorem}

\begin{proof}
Given  \eqref{volnu1},  the same proof as in Theorem \ref{lemd3}  yields
$$
 \||\nabla p_{2t}(.,y)|\|_{p}\le    \frac{C }
{\sqrt{t}\left[V(y,\sqrt{t})\right]^{1-\frac{1}{p}}},\,\ \forall\,t\ge 1,\,y\in M.
$$
On the other hand, $(R_p)_{loc}$ easily implies, by analyticity of the heat semigroup on $L^p(M,\mu)$,
 $$\||\nabla e^{- t\Delta}|\|_{p\to p} \le C\left( \frac{1}{\sqrt t}+1\right),
$$
for  all $t>0$, hence, following \cite[p.944]{ACDH},
$$
 \||\nabla p_{2t}(.,y)|\|_{p}\le    \frac{C }
{\sqrt{t}\left[V(y,\sqrt{t})\right]^{1-\frac{1}{p}}},\,\ \forall\,t\le 1,\,y\in M.
$$
One concludes as before.

\end{proof}

{\bf Remark:}
One way to ensure \eqref{volnu1} is to assume  $(D_\nu)$ and
 \begin{equation}
\sup_{x\in M}V(x,1)<+\infty.\label{noexp}
\end{equation}

\bigskip

Let us consider now the limit case $p=\infty$ in inequality  \eqref{sobpd}.

\begin{theorem}\label{lem5d}
Let $M$ satisfy $(D)$, $(DU\! E)$, and \eqref{volnu} for some $\nu>0$.
Assume
  \begin{equation}\label{sobpdi}
    \|\left|\nabla f\right|\|_{\infty}
\le  C   \|f\|_{q}^{1-\frac{\nu+q}{\alpha q}}
\|\Delta^{\alpha/2}f\|^{\frac{\nu+q}{\alpha q}} _{q},  \ \forall\,f\in
{\bf C}^\infty_0(M),
 \end{equation}
for some $q\in [1,\infty)$  and some  $\alpha > \frac{\nu}{q}+1$.
Then $(E_p)$ holds for all $p \in (1,\infty)$.
\end{theorem}

\begin{proof}
Taking again $f=p_{2t}(.,y)=\exp(-t\Delta)p_t(.,y)$, $t>0$, $y\in M$, in \eqref{sobpdi},
and using the fact that
$(U\!E)$ yields
\begin{equation}
\|p_t(.,y)\|_q\le \frac{C}{\left[V(y,\sqrt{t})\right]^{1-\frac{1}{q}}},\
\forall\,t>0,\,y\in M,\label{qq}
\end{equation}
one obtains 
$$
 \||\nabla p_{2t}(.,y)|\|_{\infty}\le    \frac{C}
{\left[V(y,\sqrt{t})\right]^{\left(1-\frac{1}{q}\right)\left(1-\frac{\nu+q}{\alpha q}\right)}} \|\Delta^{\alpha/2}\exp(-t\Delta)p_t(.,y)\|^{\frac{\nu+q}{\alpha q}} _{q},
$$
hence, by analyticity of the heat semigroup on $L^q(M,\mu)$ (when   $q=1$, it follows from $(U\!E)$,   
see for instance \cite[Lemma~9]{Da1}, \cite[Theorem~3.4.8, p.103]{Da} or \cite{Ou}),
\begin{eqnarray*}
\||\nabla p_{2t}(.,y)|\|_{\infty}
& \le&    \frac{C }
{\left[V(y,\sqrt{t})\right]^{\left(1-\frac{1}{q}\right)\left(1-\frac{\nu+q}{\alpha
q}\right)}}t^{-\frac{\alpha}{2}\left(\frac{\nu+q}{\alpha q}\right)
}\|p_t(.,y)\|^{\frac{\nu+q}{\alpha q}} _{q}\\
 &\le&    \frac{C }
{\left[V(y,\sqrt{t})\right]^{\left(1-\frac{1}{q}\right)}}t^{-\frac{\nu+q}{2
q}}\\
 &\le&    \frac{C } {\sqrt{t}V(y,\sqrt{t})}\left[t^{-\nu/2
}V(y,\sqrt{t})\right]^{1/q},\,\ \forall\,t>0,\,y\in M,
\end{eqnarray*}
hence, using \eqref{volnu},
$$\||\nabla p_{2t}(.,y)|\|_{\infty}\le  \frac{C }
{\sqrt{t}V(y,\sqrt{t})},\,\ \forall\,t>0,\,y\in M,$$
that is, $(G)$.
One concludes by applying \cite{ACDH}, namely Theorem \ref{maincor} above.

\end{proof}

\begin{theorem}\label{lem4d} Assume that $M$ has Ricci curvature bounded from below.
Let $M$ satisfy $(D)$, $(DU\! E)$, \eqref{volnu1}  and \eqref{sobpdi} for some $\nu>0$,
 some $q\in [1,\infty)$  and some  $\alpha > \frac{\nu}{q}+1$.
Then $(E_p)$ holds for all $p \in (1,\infty)$.
\end{theorem}

\begin{proof}
Given \eqref{volnu1}, the same proof as in Theorem \ref{lem5d} yields $$\||\nabla p_{2t}(.,y)|\|_{\infty}\le  \frac{C }
{\sqrt{t}V(y,\sqrt{t})},\,\ \forall\,t\ge 1,\,y\in M,$$ that is, $(G)$ for large time.
Since $M$ has Ricci curvature bounded from below, it follows from \cite{LY} that $(G)$ also holds for small time.
One concludes as before.
\end{proof}

Note that inequality \eqref{sobpdi} is  known in $\R^n$,  with $\nu=n$.

\bigskip

{\bf Remark :} 
Let the space $\mbox{Lip}(M)$ be the completion of ${\bf C}^\infty_0(M)$ with respect
to the norm
$$
\|f\|_{\mbox{Lip}}=\sup_{x\neq y}\frac{|f(x)-f(y)|}{d(x,y)}.
$$
It is well-known that, if $f\in \mbox{Lip}(M)$, then $f$ is differentiable almost everywhere and
$$
\|f\|_{\mbox{Lip}}=\|\left|\nabla f\right| \|_{\infty}.
$$

By the reiteration lemma (see for instance \cite[Proposition 2.10, p.316]{BS}),
inequality (\ref{sobpdi}) is equivalent to the embedding
$$\left[L^q_{\alpha},L^q\right]_{\theta,1}\longrightarrow 
\mbox{Lip},$$
where $\left[X,Y\right]_{\theta,r}$ denotes the real interpolation space between
$X$ and $Y$ with parameters $\theta$ and $r$,  $\theta=\frac{\nu+q}{\alpha q}$,
and $L^q_{\alpha}$ is the completion of  ${\bf C}^\infty_0(M)$ with respect to the norm
$\|\Delta^{\alpha/2}f\|_q$.
Then it is a well-known fact (see for instance \cite[Proposition~3.5.3]{BuBe} and modify it to obtain a version for homogeneous spaces)
that
$$
\left[L^q_\alpha,L^q\right]_{\theta,1}=\Lambda_{\theta\alpha}^{q,1},
$$
where  the Besov space $ \Lambda^{p,q}_\alpha$ is defined via the norm
$$
\Lambda^{q,1}_\alpha(f)=
\int_0^{+\infty}t^{k-\frac{\alpha}{ 2}}\|\Delta^ke^{-t\Delta}f\|_q\frac{dt}{t}, 
$$
for $k>\alpha/2$.
Finally \eqref{sobpdi} is equivalent to
$$\Lambda_{\frac{\nu}{ q}+1}^{q,1}\longrightarrow 
\mbox{Lip}.$$

Let us finally  consider the limit case $q=\infty$ in Theorem \ref{lem5d}.
Here, no upper volume bound is needed.

\begin{theorem}\label{lem6d}
Let $M$ satisfy $(D)$ and $(DU\! E)$. 
Assume
  \begin{equation}\label{sobpdy}
    \|\left|\nabla f\right|\|_{\infty}
\le  C   \|f\|_{\infty}^{1-\frac{1}{\alpha }}
\|\Delta^{\alpha/2}f\|^{\frac{1}{\alpha }} _{\infty},  \ \forall\,f\in
{\bf C}^\infty_0(M),
 \end{equation}
for some  $\alpha > 1$.
Then $(E_p)$ holds for all $p \in (1,\infty)$.
\end{theorem}

Let us emphasize the particular case $\alpha=2$ of inequality \eqref{sobpdy}:
\begin{equation}
\|\left|\nabla f\right|\|_{\infty}^2
\le  C   \|f\|_{\infty}
\|\Delta f\|_{\infty},  \ \forall\,f\in
{\bf C}^\infty_0(M).\label{interinf}
\end{equation}

\begin{proof}
Substituting $\exp(-t\Delta)f$ in (\ref{sobpdy}) yields
\begin{equation*}
  \|\left|\nabla \exp(-t\Delta)f\right|\|_{\infty}
\le  C   \|\exp(-t\Delta)f\|_{\infty}^{1-\frac{1}{\alpha }}
\|\Delta^{\alpha/2}\exp(-t\Delta)f\|^{\frac{1}{\alpha }} _{\infty}.
\end{equation*}
Recall that it follows from $(DU\!E)$ that  the heat semigroup is analytic on $L^1(M,\mu)$,
hence by  duality
$$\|\Delta^{\alpha/2}\exp(-t\Delta)f\| _{\infty}\le C t^{-\alpha/2 } \|f\|_{\infty}.$$
 The heat semigroup being uniformly bounded on $L^\infty(M,\mu)$, one obtains
$$  \|\left|\nabla \exp(-t\Delta)f\right|\|_{\infty}
\le  
C t^{-1/2 } \|f\|_{\infty},
$$
that is, $(G_\infty)$, or
$$\sup_{x\in M,t>0}\sqrt{t}\int_M|\nabla_xp_t(x,y)|\,d\mu(y)<\infty.$$
 It is well-known and easy to see that $(G_\infty)$ together with $(D)$ and $(U\!E)$ implies $(G)$ (in fact, these conditions are equivalent, because of the already-mentionned self-improvement of $(G)$).
One concludes again by applying \cite{ACDH}, namely Theorem \ref{maincor} above.
\end{proof}

Next we discuss a result which does not require any assumption on the volume growth and which is motivated by 
\eqref{interinf}. This result is contained in \cite{D1}, with  a similar argument, in  a discrete setting. For another approach to inequality \eqref{interinfp} below, see \cite[Section 4]{CD1}.

\begin{proposition}\label{prop12} For any $1\le p \le \infty$, 
condition $(G_p)$ is equivalent to :
\begin{equation}
\|\left|\nabla f\right|\|_{p}^2
\le  C   \|f\|_{p}
\|\Delta f\|_{p},  \ \forall\,f\in
{\bf C}^\infty_0(M).\label{interinfp}
\end{equation}
\end{proposition}
\begin{proof}
To prove that condition (\ref{interinfp}) implies $(G_p)$ we modify slightly the argument of the proof of
Theorem \ref{lem6d}. Namely,  we put $\alpha =2$ and replace $L^\infty$ norm
by $L^p$ norm.

To prove the opposite direction, write
\begin{eqnarray*}
\nabla{(I+t \Delta)^{-1}}=\int_0^\infty \nabla\exp(-s(1+t \Delta))\,ds.
\end{eqnarray*}
Hence, for suitable $f$,
\begin{eqnarray*}
\||\nabla{(I+t \Delta)^{-1}}f|\|_p\le \int_0^\infty e^{-s} \||\nabla \exp(-st \Delta)f|\|_p\,ds.
\end{eqnarray*}
Assuming $(G_p)$, one obtains, for $f\in L^q(M,\mu)$,
\begin{eqnarray*}
\||\nabla{(I+t \Delta)^{-1}}f|\|_p&\le& C \int_0^\infty e^{-s} \left(ts\right)^{-1/2} \|f\|_p \,ds \\
&=& C t^{-1/2} \|f\|_p\int_0^\infty s^{-1/2}
 e^{-s} \,ds\\
&= &C' t^{-1/2}\|f\|_p. 
\end{eqnarray*}

Hence
\begin{eqnarray*}
\||\nabla f  | \|_p&\le&
 C' t^{-1/2}\|(I+t \Delta)f\|_p \\ &\le& C'
 t^{-1/2}(\|f\|_p+ \|(t\Delta)f\|_{p}) \\ &=& C'
t^{-1/2}(\|f\|_p+t\|\Delta f\|_{p}).
\end{eqnarray*}

Taking $t=\|f\|_p\|\Delta f\|_{p}^{-1}$ yields (\ref{interinfp}).
\end{proof}

\bigskip

\section{Polynomial volume growth}\label{pg}

\begin{theorem}\label{lem4} Let $n>0$.
Suppose that $M$ satisfies  upper and lower  $n$-dimensional Gaussian
estimates
$$
ct^{-n/2}\exp\left(-\frac{d^2(x,y)}{ct}\right)\le p_t(x,y) \le Ct^{-n/2}\exp\left(-\frac{d^2(x,y)}{Ct}\right),$$
for some $C,c>0$, all $x,y \in M$ and $t>0$.

Let $p_0\in (2,\infty]$.
Then the following are equivalent:

$i)$
 \begin{equation}\label{sobp}
  \|\left|\nabla f\right|\|_p
\le  C_{p,q}
\|\Delta^{\frac{n}{2}\left(\frac{1}{q}-\frac{1}{p}\right)+\frac{1}{2}}f\|_q,\
\forall\,f\in {\bf C}^\infty_0(M),
    \end{equation}
for some $q\in (1,p)$, and all $p\in (2,p_0)$.

$ii)$
  $(R_p)$ holds for all $p \in (2,p_0)$.
\end{theorem}

\begin{proof} Let $q$ and $p$ be such that $1<q<p<\infty$ and $p>2$.
According to \cite{V}, the following  Sobolev inequality is a   consequence of the upper heat kernel estimate :
$$ \| f\|_p
\le  C \|\Delta^{\frac{n}{2}\left(\frac{1}{q}-\frac{1}{p}\right)}f\|_q,\
\forall\,f\in {\bf C}^\infty_0(M),
$$
and in particular
$$ \| \Delta^{1/2}f\|_p
\le  C
\|\Delta^{\frac{n}{2}\left(\frac{1}{q}-\frac{1}{p}\right)+\frac{1}{2}}f\|_q,\
\forall\,f\in {\bf C}^\infty_0(M).
$$
Thus $(R_p)$ for some $p>2$ implies \eqref{sobp} for all $q$ such that
$1<q<p$,
and in particular $ii)$ implies $i)$.

Conversely,  observe that the heat kernel estimates imply $V(x,r) \simeq r^n$, $\forall\,r>0,\,x\in M$
(see \cite[Theorem 3.2]{G}).
Therefore Theorem \ref{lemd3} applies with $\nu=n$ and shows that
$i)$ implies $ii)$.

\end{proof}

Remarks similar to those after Theorem \ref{lemd3} are in order.  We add one more.

\bigskip

{\bf Remark:} 
According to \cite[Theorem~0.4]{AC}, under the assumptions of Theorem \ref{lem4}, there  always exists a $p_0$
such that $ii)$ holds. It would be nice to have a proof of this fact using $i)$.
 
 \bigskip

Again, we shall now consider the limit case $p=\infty$ of inequality \eqref{sobp}.

We shall have to make local assumptions in order to ensure that the quantity
$$
\theta(t):=\sup_{0<u\le t,\,x\in M} u^{n/2}p_{u}(x,x)=\sup_{0<u\le t} u^{n/2}\|\exp(-u\Delta)\|_{1\to \infty}
$$
 is finite for some (all) $t>0$. For instance, a local Sobolev inequality of dimension $n$ is enough, since then $\sup_{\in M}p_t(x,x)\le Ct^{-n/2}, \ 0<t\le 1$.
This holds for instance if $\dim M \le n$, $M$ has Ricci curvature bounded from below  and satisfies the matching condition to \eqref{noexp} : $\inf_{x\in M}V(x,1)>0$.

\begin{theorem}\label{lem3} Assume that $M$ has Ricci curvature bounded from below.
Let $n\in \N^*$. Assume that \begin{equation}
V(x,r) \simeq r^n,\ \forall\,r>0,\,x\in M.\label{voln}
\end{equation}
Then
$M$ satisfies the heat kernel gradient estimate $(G)$,
that is  \begin{equation}
\left|\nabla_xp_t(x,y)\right| \le Ct^{-\frac{n+1}{2}}\exp\left(-\frac{d^2(x,y)}{Ct}\right),\label{ggn}
\end{equation}
for some $C>0$, all $x,y \in M$ and $t>0$,
if and only if
 \begin{equation}\label{1aq}
  \|\left|\nabla f\right|\|_{\infty}
\le  C   \|f\|_{q}^{1-\frac{n+q}{\alpha q}}
\|\Delta^{\alpha/2}f\|^{\frac{n+q}{\alpha q}} _{q},  \ \forall\,f\in {\bf
C}^\infty_0(M),
    \end{equation}
    for some (all) $q\in (1,\infty)$ and some (all) $\alpha > \frac{n}{q}+1$.
   Moreover if \eqref{ggn} or \eqref{1aq} holds, then $M$
 satisfies $(LY)$, that is the upper and lower  $n$-dimensional Gaussian
estimates
\begin{equation}
ct^{-n/2}\exp\left(-C\frac{d^2(x,y)}{t}\right)\le p_t(x,y) \le Ct^{-n/2}\exp\left(-c\frac{d^2(x,y)}{t}\right),\label{gn}
\end{equation}
for some $C,c>0$, all $x,y \in M$ and $t>0$.
\end{theorem}

\bigskip

The following result is a direct consequence of Theorem \ref{lem3} together with  \cite[Theorem 1.4]{ACDH}, that is, Theorem \ref{maincor} above.

\begin{coro}\label{E} Assume that $M$ has Ricci curvature bounded from below,  and satisfies \eqref{voln} and \eqref{1aq}.
Then $(E_p)$ holds
for all $p \in (1,\infty)$.
\end{coro}

Let us prepare the proof of Theorem \ref{lem3} with two   lemmas. The first one is reminiscent of Proposition \ref{gr} : it shows that
a certain gradient estimate implies an upper bound of the heat kernel.

\begin{lemma}\label{grp} Assume that $M$ has Ricci curvature bounded from below  and let $n>0$.
Assume that, for some $c>0$,
\begin{equation}
V(x,r)\ge cr^n,\label{vl}
\end{equation}
for all $x\in M$ and $r>0$. Next
suppose that, for some $q\in[1,\infty]$,
$$
\||\nabla \exp(-t\Delta)|\|_{q\to \infty} \le C
t^{-\frac{n+q}{2q}}\leqno{(G_{q,\infty}^n)}
$$
for all $t>0$.
Then there exists a constant $C'$ such that
$$
\sup_{x\in M}p_t(x,x)  \le  C't^{-n/2}
$$
for all $t>0$.
\end{lemma}
\begin{proof}
Set
$
\theta(t)=\sup_{0<u\le t,\,x\in M} u^{n/2}p_{u}(x,x)=\sup_{0<u\le t} u^{n/2}\|\exp(-u\Delta)\|_{1\to \infty}.
$

Remember that  the curvature assumption   together with the volume lower bound ensures the finiteness of $\theta(t)$ for all $t>0$.
Using $(G_{q,\infty}^n)$ and  interpolation, we can write
\begin{eqnarray*}
\||\nabla \exp(-2s\Delta)|\|_{1\to \infty} & \le& \||\nabla
\exp(-s\Delta)|\|_{q\to \infty} \|\exp(-s\Delta)\|_{1 \to q} \\
& \le & C s^{-\frac{n+q}{2q}}\left(\theta(s) s^{-\frac{n}{2}}\right)^{1-\frac{1}{q}}
\\
&=& C s^{-\frac{n+1}{2}}\theta(s)^{1-\frac{1}{q}}.
\end{eqnarray*}
For $x\in M$ and $s>0$, define
$$
K=K(s,x)=\frac{s^{n/2}p_{2s}(x,x)}{2}.
$$
For all $y\in B\left(x,\frac{K\sqrt{s}}{C\theta(s)^{1-\frac{1}{q}}}\right)$,
\begin{eqnarray*}
|p_{2s}(y,x)-p_{2s}(x,x)| &\le& {d(y,x)} \sup_{z\in M}|\nabla p_{2s}(z,x)|\\&\le& {d(y,x)} \||\nabla \exp(-2s\Delta)|\|_{1 \to
\infty}\\ & \le&
\frac{K\sqrt{s}}{C\theta(s)^{1-\frac{1}{q}}}\frac{C\theta(s)^{1-\frac{1}{q}}}{s^{\frac{n+1}{2}}}=\frac{K}{s^{n/2}
},
\end{eqnarray*}
therefore
$$
p_{2s}(y,x) \ge p_{2s}(x,x)-\frac{K}{s^{n/2} }=\frac{2K}{s^{n/2} } -\frac{K}{s^{n/2} }=\frac{K}{s^{n/2} }.
$$
Hence
\begin{eqnarray*}
 1&\ge&\int_M p_{2s}(y,x)\,d\mu(y)   \ge
\int_{B\left(x,\frac{K\sqrt{s}}{C\theta(s)^{1-\frac{1}{q}}}\right)} p_{2s}(y,x)\,d\mu(y)\\&\ge&
\int_{B\left(x,\frac{K\sqrt{s}}{C\theta(s)^{1-\frac{1}{q}}}\right)} \frac{K}{s^{n/2}} \,d\mu(y) \ge
c  \frac{K^{n+1}}{C^n \theta(s)^{n(1-\frac{1}{q})} },
 \end{eqnarray*}
using (\ref{vl}) in the last inequality.
Since $\theta$ is obviously non-decreasing, we also
 have
 $$
 1 \ge
c  \frac{K^{n+1}}{C^n \theta(2s)^{n(1-\frac{1}{q})}  },
 $$
that is
$$K(s,x)\le \left( \frac{C^n}{c}\theta(2s)^{n(1-\frac{1}{q})}
\right)^{\frac{1}{n+1}},$$
hence
$$K(s,x)\le \left( \frac{C^n}{c}\theta(2t)^{n(1-\frac{1}{q})}
\right)^{\frac{1}{n+1}}$$
for $0<s\le t$.

Taking supremum in $x$ and $s$ yields
$$2^{-\frac{n}{2}-1}\theta(2t)\le \left( \frac{C^n}{c}\theta(2t)^{n(1-\frac{1}{q})}
\right)^{\frac{1}{n+1}}.$$

Since $\frac{n(1-\frac{1}{q})}{n+1}<1$,
it follows that $\theta$ is bounded from above, which proves the claim.
\end{proof}

{\bf Remark:} One can write a version of the above lemma in the case where 
$V(x,r)\ge v(r)$, for some doubling function $v$.

\bigskip

The lemma below yields as a by-product a new  proof of inequality \eqref{1aq} in $\R^n$. It does not require any volume growth assumption.

\begin{lemma}\label{lem1}
Let $1 < q < \infty $ and $n>0$. The 
following estimates are equivalent:

i) $$
\||\nabla \exp(-t\Delta)|\|_{q\to \infty} \le C
t^{-\frac{n+q}{2q}},\ \forall\,t>0.\leqno{(G_{q,\infty}^n)}
$$

ii)    $$
  \||\nabla   (I+t \Delta)^{-\alpha/2}  |\|_{q\to \infty}
 \le    C_\alpha t^{-\frac{n+q}{2q}} ,$$
    for some (all) $\alpha > \frac{n}{q}+1$ and all $t>0$.
    
iii)    \begin{equation*}
  \|\left|\nabla f\right|\|_{\infty}
\le  C   \|f\|_{q}^{1-\frac{n+q}{\alpha q}}
\|\Delta^{\alpha/2}f\|^{\frac{n+q}{\alpha q}} _{q},  \ \forall\,f\in {\bf
C}^\infty_0(M),
    \end{equation*} for  some (all) $\alpha > \frac{n}{q}+1$, that is,  \eqref{1aq}.
\end{lemma}
\begin{proof}
We shall show that $i)\Rightarrow ii)\Rightarrow iii)\Rightarrow i)$.

Write
\begin{eqnarray*}
\nabla{(I+t \Delta)^{-\alpha/2}}=\int_0^\infty s^{(\alpha/2) -1
}\nabla\exp(-s(1+t \Delta))\,ds.
\end{eqnarray*}
Hence, for suitable $f$,
\begin{eqnarray*}
\||\nabla{(I+t \Delta)^{-\alpha/2}}f|\|_\infty\le \int_0^\infty s^{(\alpha/2) -1
}e^{-s} \||\nabla \exp(-st \Delta)f|\|_\infty\,ds.
\end{eqnarray*}
Assuming $(G_{q,\infty}^n)$, one obtains, for $f\in L^q(M,\mu)$,
\begin{eqnarray*}
\||\nabla{(I+t \Delta)^{-\alpha/2}}f|\|_\infty&\le& C \int_0^\infty s^{(\alpha/2)
-1 }e^{-s} \left(ts\right)^{-\frac{n+q}{2q}} \|f\|_q \,ds \\
&=& C t^{-\frac{n+q}{2q}} \|f\|_q \int_0^\infty s^{(\alpha/2) -\frac{n}{2q} -(3/2)}
 e^{-s} \,ds\\
&= &C_\alpha t^{-\frac{n+q}{2q}}\|f\|_q, 
\end{eqnarray*}
since $\alpha>\frac{n}{q}+1$.

Assume $ii)$, and write
\begin{eqnarray*}
\||\nabla f| \|_{\infty}&\le&
 C t^{-\frac{n+q}{2q}}\|(I+t \Delta)^{\alpha/2}f\|_q \\ &\le& C_\alpha
 t^{-\frac{n+q}{2q}}(\|f\|_q+ \|(t\Delta)^{\alpha/2}f\|_{q}) \\ &=& C
t^{-\frac{n+q}{2q}}(\|f\|_q+t^{\alpha/2}\|\Delta^{\alpha/2}f\|_{q}).
\end{eqnarray*}
The second inequality relies on the $L^p$-boundedness of the operator
$(I+t \Delta)^{\alpha/2}(I+(t \Delta)^{\alpha/2})^{-1}$ (see \cite{S}, or use analyticity).

Taking $t=\|f\|_q^{2/\alpha}\|\Delta^{\alpha/2}f\|_{q}^{-2/\alpha}$ yields
$iii)$.

Finally, assume $iii)$.
Replacing  $f$ by $\exp(-t \Delta)f$, one obtains, by  contractivity and
analyticity of the heat semigroup on $L^q(M,\mu)$,
\begin{eqnarray*}
\||\nabla{\exp(-t \Delta)}f|\|_{\infty}&\le&C \|{\exp(-t
\Delta)}f\|_{q}^{1-\frac{n+q}{\alpha q}}
\|\Delta^{\alpha/2}{\exp(-t \Delta)}f\|_{q}^{\frac{n+q}{\alpha q}}\\
&\le& C t^{-\alpha\frac{n+q}{2\alpha q}}\|f\|_q=Ct^{-\frac{n+q}{2q}}\|f\|_q,
\end{eqnarray*}
that is, $i)$.  
\end{proof}

{\bf Remark:}  
The above lemma also holds for $q=1,\infty$, provided the heat semigroup is analytic on $L^1(M,\mu)$, which is the case, as we already said, if it satisfies Gaussian estimates and $(D)$ holds.  Note that, according to Lemma \ref{grp}, this is automatic from $i)$ under the boundedness from below of the Ricci curvature and  \eqref{voln}.

\begin{proof}
[Proof of  Theorem~\ref{lem3}] Assume (\ref{1aq}).
By Lemma~\ref{lem1},  $(G_{q,\infty}^n)$ follows, and by Lemma \ref{grp},
\begin{equation}
\sup_{x\in M}p_t(x,x)= \|\exp(-t\Delta)\|_{1\to \infty} \le  C't^{-n/2}\label{uq}
\end{equation}
for all $t>0$. The  Gaussian upper  bound follows:
$$p_t(x,y) \le Ct^{-n/2}\exp\left(-\frac{d^2(x,y)}{Ct}\right),$$
for some $C,c>0$, all $x,y \in M$ and $t>0$.
By interpolation, \eqref{uq} yields
\begin{equation}
 \|{\exp(-t
\Delta)}\|_{1\to q}\le Ct^{-n(1-\frac{1}{q})/2}.\label{uqq}
\end{equation}
Combining $(G_{q,\infty}^n)$ with \eqref{uqq} yields
\begin{equation}
\sup_{x,y\in M} |\nabla_x p_t(x,y)|=\||\nabla \exp(-t\Delta)|\|_{1\to \infty} \le C
t^{-\frac{n+1}{2}},\label{ggo}
\end{equation}
which  together with the upper  bound yields the Gaussian lower bound
$$
ct^{-n/2}\exp\left(-\frac{d^2(x,y)}{ct}\right)\le p_t(x,y).$$
Finally,   \eqref{ggo} self-improves into \eqref{ggn}.

Conversely,   \eqref{ggn} obviously  implies  $(G_{1,\infty}^n)$ and, together with the volume upper bound,
$(G_{\infty,\infty}^n)$ or, in other words, $(G_\infty)$:
$$
\sup_{x\in M}\int_M|\nabla_xp_t(x,y)|\,d\mu(y)=\||\nabla \exp(-t\Delta)|\|_{\infty\to \infty} \le C
t^{-1/2}.
$$
By interpolation, one obtains $(G_{q,\infty}^n)$, therefore \eqref{1aq}, thanks to Lemma \ref{lem1}.
\end{proof}

{\bf Remarks:} 

-As a consequence of Theorem \ref{lem3}, \eqref{1aq} implies, using the results in \cite{CGN},
$$ \| f\|_{\infty}
\le  C   \|f\|_{q}^{1-\frac{n}{\alpha q}}
\|\Delta^{\alpha/2}f\|^{\frac{n}{\alpha q}} _{q},  \ \forall\,f\in {\bf
C}^\infty_0(M),
$$
for $\alpha>n/q$, $q\in [1,\infty)$,
 and, using the results in \cite{C},
$$\frac{|f(x)-f(y)|}{\left[d(x,y)\right]^{\alpha-(n/q)}}\le C \|\Delta^{\alpha/2}f\|_{q},  \ \forall\,x,y\in M, \,f\in {\bf
C}^\infty_0(M),$$
for $\alpha>n/q$, $q\in [1,\infty)$.
It would be interesting to have a direct  proof of these two implications.

-According to known results on Riesz transforms (see \cite{ACDH} for references), \eqref{1aq}  is true for manifolds with non-negative Ricci curvature, Lie groups with polynomial volume growth, cocompact coverings with polynomial volume growth. Again, it would be interesting to have direct proofs.

-It would interesting to study the stability under perturbation of inequalities \eqref{1aq} or \eqref{sobp}, in the light of the result in \cite{CN}.  

\section{Applications}\label{ap}
Now we consider  a uniformly elliptic operator $H$ in divergence form acting on $\R^n$,  $n\in \N^*$, that is
$$
Hf=-\sum_{i,j=1}^n\partial_i(a_{ij}\partial_jf)
$$  
where $a_{ij}\in L^\infty$ for all $1\le i,j \le n$, and the matrix $(a_{ij}(x))_{1\le i,j \le n}$  is a symmetric matrix with real coefficients, such that
$$
\sum_{i,j}a_{ij}(x)\xi_j\xi_i \ge c|\xi|^2, \mbox{ for a.e. } x,\xi\in\R^n,
$$
for some $c>0$.
Next  let $\Delta$  denote the standard non-negative Laplace operator acting on $\R^n$.

It follows from the above uniform ellipticity 
assumption and the boundedness of the coefficients that 
$$
|\nabla_H f(x)|^2=\sum_{i,j}a_{ij}(x)\partial_jf(x)\partial_i f(x) \simeq  |\nabla f(x)|^2.
$$
We say that   $H$ satisfies $(R_p)$  for some $p\in (1,\infty))$ if
$$
  \||\nabla_H f|\|_{p}
\le  C_p
 \|H^{1/2}f\|_{p}, \ \forall f\in{\bf C}^\infty_0(\R^n),
$$
which according to the above remark is equivalent
to
$$
  \|\left|\nabla f\right|\|_{p}
\le  C_p
 \|H^{1/2}f\|_{p}, \ \forall f\in{\bf C}^\infty_0(\R^n).
$$
To avoid technicalities  we assume in what follows that all coefficients  $a_{ij}, b_{ij}$ discussed below
are smooth. However we point out that this assumption  
can be substantially  relaxed.

Recall that  the Gaussian estimates do hold for $e^{-tH}$ and that the above framework applies.

The  assumption in our first application may be seen as some boundedness  for  the higher order Riesz transform associated with $H$.

\begin{theorem}\label{appl1}
Suppose that 
\begin{equation}\label{high1}
\|\Delta^{\alpha/2} f\|_{q_0}\le C \|H^{\alpha/2} f\|_{q_0},\ \forall f\in{\bf C}^\infty_0(\R^n),
\end{equation}
for some $\alpha >1$ and $1<q_0<\infty$. Then, if  $\alpha <  \frac{n}{q_0}+1$, $H$ satisfies $(R_p)$   
for all $p\in(1, p_0)$, where $p_0=\frac{n}{\frac{n}{q_0}+1-\alpha}$,
and  if $\alpha \ge \frac{n}{q_0}+1$, 
 $H$ satisfies $(R_p)$  for all $p\in (1,\infty)$.
\end{theorem}
\begin{proof}
The  boundedness of  the classical Riesz transform on $L^{p}(\R^n,dx)$ together with the Sobolev inequality in $\R^n$  imply, for $1<q_0<p<\infty$, 
$$
   \|\left|\nabla f\right|\|_{p} \le  C
\|\Delta^{\frac{n}{2}\left(\frac{1}{q_0}-\frac{1}{p}\right)+\frac{1}{2}}f\|_{q_0}
\le  C'
\|\Delta^{\alpha/2}f\|_{q_0}^\theta\|f\|_{q_0}^{1-\theta}, \ \forall f\in{\bf C}^\infty_0(\R^n),
$$
as soon as $\alpha\ge  n(\frac{1}{q_0}-\frac{1}{p})+1$, $\theta\in (0,1]$ being such that
$\alpha\theta=  n(\frac{1}{q_0}-\frac{1}{p})+1$.
Now let $\alpha$ and $q_0$ be such that  \eqref{high1} holds.
If $\alpha \ge \frac{n}{q_0}+1$,  choose any $p>q_0$.
If $\alpha < \frac{n}{q_0}+1$,  choose $p_0$  in $(q_0,\infty)$ so that  $\alpha=n(\frac{1}{q_0}-\frac{1}{p_0})+1$.
In both cases,
$$
   \|\left|\nabla f\right|\|_{p_0} \le    C
\|H^{\alpha/2}f\|_{q_0}^\theta\|f\|_{q_0}^{1-\theta}, \ \forall f\in{\bf C}^\infty_0(\R^n),
$$
and $(R_p)$ for $1<p<p_0$ follows from Theorem~\ref{lemd3} and the  remarks afterwards. 
\end{proof}

Our next application says that  $(R_p)$ also holds for  small $L^\infty \cap W^{1,n}$ perturbations of  operators with bounded second order Riesz transform.

To state this result   we set
$$
H_\varepsilon f =H f+\varepsilon \sum_{i,j}\partial_ib_{ij}(x)\partial_jf,
$$
where $H=H_0$ is as above.
We do not assume here that the matrix $(b_{ij}(x))_{1\le i,j \le n}$  is positive definite. However we assume
that  $b_{ij}\in L^\infty(\R^n,dx)$ and that $\varepsilon$ is small enough so that the operator $H_\varepsilon$ is uniformly elliptic. 
\begin{theorem}\label{appl2}
Suppose that $b_{ij}\in L^\infty(\R^n,dx)$ for all $1\le i,j \le n$ and that $\partial_jb_{ij}\in L^n(\R^n,dx)$ for all $1\le i,j \le n$.
Next assume that for some $q_0< n$
\begin{equation}\label{second}
\|\Delta f\|_{q_0}\le C_{q_0} \|H_0 f\|_{q_0},\ \forall f\in{\bf C}^\infty_0(\R^n),
\end{equation}
Then there exists $\gamma>0$ such that  {\rm (\ref{second})} holds for $H_\varepsilon$ for $\varepsilon <\gamma$ and so
$(R_p)$ holds for $H_\varepsilon$ for all $\varepsilon <\gamma$ and 
 $1<p<p_0$ where $\frac{1}{p_0}+\frac{1}{n}=\frac{1}{q_0}$.  
\end{theorem}
\begin{proof}
Note that (\ref{second}) is just condition (\ref{high1}) for $\alpha=2$. 
We are going to prove that this inequality extends from $H$ to $H_\varepsilon$ for $0<\varepsilon<\gamma$ and apply
Theorem~\ref{appl1}. To this purpose, it is enough to show that for some $\gamma>0$ and for all $\varepsilon <\gamma$
\begin{equation}\label{dif}
\|H_\varepsilon f -H_0f\|_{p_0} \le \frac{1}{2C_{q_0}} \|\Delta f\|_{p_0}.
\end{equation}

Now
$$
\|H_\varepsilon f -H_0 f\|_{q_0}\le \varepsilon \left( \sum_{i,j}\| \partial_ib_{ij}\partial_jf\|_{q_0}   \right).
$$
Since 
$$ \partial_i(b_{ij}\partial_jf)=b_{ij}(\partial_i\partial_jf)+(\partial_ib_{ij})(\partial_jf),$$
one may write
$$\sum_{i,j}\| \partial_ib_{ij}\partial_jf\|_{q_0} \le 
\max_{i,j}\|b_{ij}\|_\infty\sum_{i,j}\| \partial_i\partial_jf\|_{q_0}+\sum_{i,j}\| (\partial_ib_{ij})(\partial_jf)\|_{q_0},
$$
hence
\begin{eqnarray*}
&&\sum_{i,j}\| \partial_ib_{ij}\partial_jf\|_{q_0} \\
&\le&\max_{i,j}\|b_{ij}\|_\infty\sum_{i,j}\| \partial_i\partial_jf\|_{q_0}+\sum_{i,j}\| (\partial_ib_{ij})(\partial_jf)\|_{q_0}\\ &\le& n^2\max_{i,j}\|b_{ij}\|_\infty \| \partial_i\partial_j\Delta^{-1} \|_{q_0 \to q_0}\|\Delta f\|_{q_0}
+\sum_{i,j}\| \partial_ib_{ij}\|_n\|\partial_j f\|_{p_0}.
\end{eqnarray*}

Here we have used the $L^{q_0}$ boundedness of the second order Riesz transform in $\R^n$
and the H\"older inequality  $\|fg\|_{q_0} \le \|f\|_n\|g\|_{p_0}$.

Now recall that an inequality similar to (\ref{sobp})  holds in $\R^n$, that is $$\|\partial_jf\|_{p_0}\le C\| \Delta f \|_{q_0}.$$ 

Therefore
\begin{eqnarray*}&&\sum_{i,j}\| \partial_ib_{ij}\partial_jf\|_{q_0}\le \\&& \left(n^2\max_{i,j}\|b_{ij}\|_\infty \| \partial_i\partial_j\Delta^{-1} \|_{q_0 \to q_0}   +\max_i
\| \partial_i\Delta^{-1} \|_{q_0 \to p_0} \sum_{i,j}\| \partial_ib_{ij}\|_n   \right) 
\|\Delta f\|_{q_0}.   
\end{eqnarray*}
This yields (\ref{dif}) with 
$$\gamma \left(n^2\max_{i,j}\|b_{ij}\|_\infty \| \partial_i\partial_j\Delta^{-1} \|_{q_0 \to q_0}  + 
\max_i\| \partial_i\Delta^{-1} \|_{q_0 \to p_0} \sum_{i,j}\| \partial_ib_{ij}\|_n   \right) =\frac{1}{2C_{q_0}}.$$
\end{proof}

\bigskip

The second order Riesz transform bound  (\ref{second}) is known  
for various
large classes of operators. We discuss one instance of such class next.

{\bf Example:}
\it  Assume the coefficients $a_{ij}$ of $H$ are  continuous and
periodic with a common period and that $\sum^n_{i=1}\partial_ia_{ij}=0$ for $1\le j \le n$. Then
$$
\|Hf\|_p\simeq \|\Delta f\|_{p}, \ \forall f\in{\bf C}^\infty_0(\R^n),
$$
for all $1 \le p \le \infty$ (see  \cite[Theorem 1.3]{ERS}),  so that $H$ satisfies the assumption of Theorem~{\rm \ref{appl2}}. In \cite{ERS}, it is proved that $(R_p)$ holds for such $H$, but the above shows that is it also holds for  small $L^\infty \cap W^{1,n}$ perturbations of  $H$.
\rm

\bigskip

{\bf Remark:}
It is interesting to compare Theorem~\ref{appl1},  which proves that boundedness of second 
order  Riesz transform implies boundedness of first
order Riesz  transform on a larger range on $L^p$ spaces,  with the results obtained in \cite{MM}.  See also \cite[(1.26)]{HMM}.

\bigskip

{\bf Acknowledgements: } The authors would like to thank Professor Vladimir Maz'ya for interesting discussions about the proof of inequality  \eqref{1aq} in the Euclidean space.

\end{document}